
\input amstex

\overfullrule=0pt

\documentstyle{amsppt}

\topmatter

\title\nofrills Bounding canonical genus bounds volume \endtitle
\author  Mark Brittenham \endauthor

\leftheadtext\nofrills{Mark Brittenham}
\rightheadtext\nofrills{Bounding canonical genus bounds volume}

\affil   University of  North Texas \endaffil
\address   Department of Mathematics, University of North Texas, Denton, TX 76203  \endaddress
\email   britten\@unt.edu \endemail
\thanks   Research supported in part by NSF grant \# DMS$-$9704811 \endthanks

\abstract In this paper we show that there is an upper bound on 
the volume of a
hyperbolic knot in the 3-sphere with canonical genus $g$. 
This bound
can in fact be chosen to be linear in $g$. In other 
words, if Seifert's algorithm builds a surface with small genus 
for a hyperbolic knot, then the complement of the knot cannot have large 
volume. \endabstract

\endtopmatter

\document

\heading{\S 0 \\ Introduction}\endheading

A Seifert surface for a knot $K$ in the 3-sphere is an embedded
 orientable surface $\Sigma$, whose boundary equals $K$.
In 1934, Seifert [Se] gave a very simple algorithm for 
constructing a Seifert surface for a knot $K$, from a 
diagram, or projection, $D$ of the knot. Thus every knot 
has a Seifert surface.

This fact gave rise to one of the
first invariants for knots; the {\it genus} of a knot $K$, $g(K)$, is the 
minimum among all genera of Seifert surfaces for $K$. It is trivially
an invariant of the knot, since we consider all (oriented)
surfaces with boundary a knot isotopic to $K$. By ambient isotopy, this 
is equivalent to considering all surfaces with boundary equal to $K$.
This invariant is, however, remarkably hard to compute, in
part because it in principal requires us to construct all of the
Seifert surfaces for the knot, in order to be certain that 
we have built its least genus representative. In some circumstances,
however, one can show that Seifert's algorithm itself will build a 
least genus Seifert surface, e.g., starting with a reduced alternating
projection of an alternating knot [Cr],[Mu].

In general, however, it seems that Seifert's algorithm rarely builds
a least genus Seifert surface for a knot. This leads us to to define a 
new invariant; the {\it canonical genus}, $g_c(K)$, of a knot $K$ 
is the minimum
among all genera of Seifert surfaces built by Seifert's algorithm
on a diagram of the knot $K$. Since the minimum is taken over fewer
surfaces (not all Seifert surfaces, as we shall see, can be built by
Seifert's algorithm), we immediately have $g(K)\leq g_c(K)$ for all
knots $K$.

There is a third notion of genus, 
related to the previous two, which is based of the fact that
every surface $\Sigma$ built by Seifert's algorithm has the property 
that $S^3\setminus\Sigma$ is a handlebody, i.e. has free fundamental
group. Such a Seifert surface is called {\it free}, and we can 
then define the {\it free genus}, $g_f(K)$, of a knot $K$ as the 
minimum among all genera of free Seifert surfaces for $K$. The
considerations above then show that for any knot $K$, 
$g(K)\leq g_f(K) \leq g_c(K)$.

It has been previously shown that each of these genera are
distinct, i.e., there are knots for which these genera differ. 
Morton [Mn] showed that there are knots whose 
genus and canonical genus differ, using a relationship between
the canonical genus and the degree of the Jones polynomial of a knot.
Moriah [Mh] showed that for most doubled knots $K$, $g(K)$=1
but $g_f(K)$ is large, using a relationship between $g_f(K)$ and the 
minimum number of generators for the fundamental group of the 
exterior of $K$. More recently, Kobayashi and Kobayashi [KK] have
shown that for $K$ a connected sum of $n$ doubles of the trefoil
knot, $g(K)=n$, $g_f(K)=2n$, and $g_c(K)=3n$. 

Our interest in these notions of genus was prompted by the 
following question: if a $\Sigma$ is an {\it incompressible} Seifert
surface for a knot $K$, whose complement is a handlebody, can it
be constructed by applying Seifert's algorithm to some projection of the
knot $K$? On the face of it, the answer should be `No', since the
surface $\Sigma$ might be a free genus-minimizing surface; if the 
knot had higher canonical genus, then $\Sigma$ could never be built 
via Seifert's algorithm. Unfortunately, for the only examples known
where $g_f(K)$ and $g_c(K)$ differ [KK], the free genus-minimizing
surfaces are all compressible! In looking for examples to settle this
question, the key ingredient that we found is a relationship 
between the canonical genus of a (hyperbolic) knot and the volume of
its complement.


\proclaim{Theorem} For any $g$, there is a finite collection of
hyperbolic
links $L_1$,...,$L_k$ in the 3-sphere $S^3$, so that for any 
hyperbolic knot $K$ in $S^3$ with canonical genus less than or 
equal to $g$, $K$ can be obtained 
by $1/n_i$ Dehn surgeries on the unknotted components of 
one of the $L_i$.\endproclaim

By a result of Thurston [Th] (see also [Ag]), the hyperbolic volume 
of the knot $K$ 
must then be less than the volume of the corresponding link $L_i$. 
This will allow us to show:

\proclaim{Corollary} If $K$ is a hyperbolic knot with canonical genus $g$,
then the hyperbolic volume of the complement of $K$, vol$X(K)$, is 
less than $120gV_0$, where $V_0$ is the volume of the hyperbolic regular ideal 
tetrahedron.\endproclaim

This result constitutes half of a program to find hyperbolic 
knots whose canonical genus differs from their free genus; 
in particular, to find knots with very low free genus but very high 
canonical genus. In a sequel to this paper [Br], we find a family
of hyperbolic knots with free genus one, but arbitrarily large 
volume. The results of this paper show that the canonical 
genus of these knots must then be arbitrarily large, as well. 
But a free genus one Seifert surface for a 
non-trivial knot must be incompressible. Otherwise, 
compressing it gives a genus zero surface, i.e., a disk, so the 
knot will be trivial.
So these knots also provide examples of knots with incompressible
free Seifert surfaces which cannot be obtained from Seifert's algorithm.

\heading{\S 1 \\ Outline and preliminary maneuvers}\endheading

Our basic approach will be to start with a hyperbolic knot having 
canonical genus $g$, and, by changing crossings, replace it with
a hyperbolic \underbar{alternating} knot with canonical genus $g$.
We will then replace the alternating knot with an augmented 
alternating link in the sense of [Ad], that is, a link consisting 
of an underlying alternating knot together with unknotted loops
going around pairs of arcs near some of the crossings (see Figure 1). We shall 
see that the original hyperbolic knot can be retrieved from this link 
by doing $1/n_i$ Dehn surgeries on each of the unknotted components.
We shall also see that this replacement process will yield only 
finitely many augmented alternating links. This is the collection 
of links claimed in the theorem.

\input epsf.tex

\leavevmode

\epsfxsize=2in
\centerline{{\epsfbox{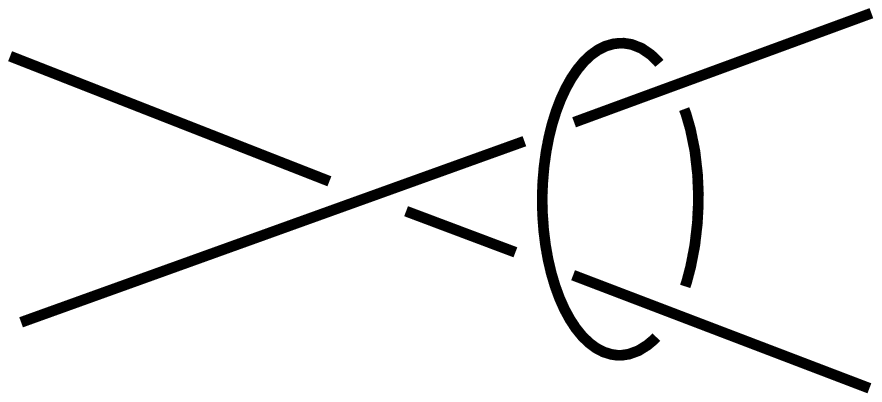}}}

\centerline{Figure 1}

We now start with a knot $K$ with canonical genus $g$, and let $D$ be a diagram 
for $K$ for which 
Seifert's algorithm builds a surface $\Sigma$ of genus $g$. We briefly recall
the outline of Seifert's algorithm. Choosing an orientation on $K$,
we remove each crossing in $D$ and glue the resulting four ends together 
according to the orientation of $K$ (see Figure 2), obtaining a collection of disjoint 
Seifert circles. The circles inherit an orientation from the
orientation of $K$. These circles may be nested, but by imagining the one lying 
inside of another to be slightly higher,
the circles then bound disjoint disks, each inheriting a (normal) orientation
from the orientation of its boundary. These can be connected by half-twisted
bands as dictated by the original crossings, to obtain a Seifert surface for 
our original knot; see Figure 2. The resulting surface is normally oriented, hence
orientable, because the half-twisted bands always join (disjoint) disks of opposite
normal orientation, or (nested) disks of the same orientation.

\leavevmode

\epsfxsize=4.9in
\centerline{{\epsfbox{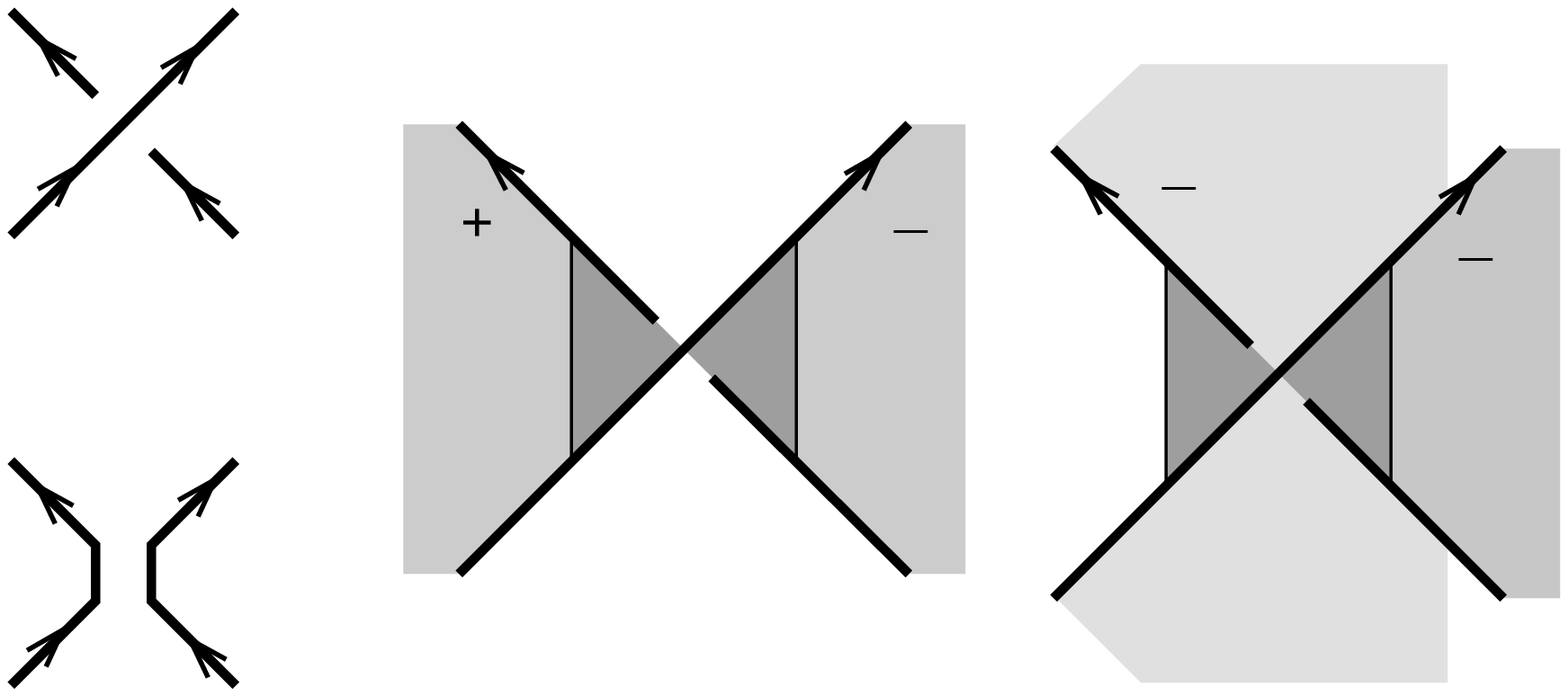}}}

\centerline{Figure 2}


In building our collection of links,
the first replacement consists of changing the crossings of our
diagram $D$ to those of an alternating diagram $D^\prime$, that is,
the diagram $D^\prime$ of an alternating knot $K^\prime$. (The proof that
this can be done is an interesting exercise in combinatorics; 
essentially one needs
to note that, starting at an arc of a crossing and walking along the knot, 
the number of crossings
we pass through before first returning to our starting point is always even. 
This, in turn, can be seen by induction, since a shortest such path traces out 
a loop in the projection plane, bounding a disk, which
meets the knot in arcs, having an even number of endpoints.)
Changing the crossings of a diagram does not change the genus of 
the surface that Seifert's algorithm builds, merely the direction in which 
some of the half-twisted bands are twisted.

Our only problem is that the resulting knot $K^\prime$ might not 
be hyperbolic. But by [Me], $K^\prime$ can fail to be hyperbolic in only
two ways. First, $D^\prime$ might be a (2,$k$)-torus knot diagram; but then so is 
$D$ (for a different $k$). Second, $K^\prime$ might be 
representable as a connected 
sum of two knots $K_1$ and $K_2$, by a loop $\gamma$ in the projection 
plane crossing the diagram 
$D^\prime$ exactly twice. But since $K$ is prime, the diagram $D$ cut along $\gamma$ 
gives a knot isotopic to $K$ and the unknot (see Figure 3). 
Seifert's algorithm applied to $D^\prime$ is a boundary connected sum of canonical Seifert
surfaces for $K_1$ and $K_2$ (which are least genus), while applied to $D$ it gives a 
boundary connected sum of
canonical Seifert surfaces for $K$ and the unknot. Consequently, $K_1$ is obtained
from a diagram of $K$ by changing crossings, and Seifert's algorithm gives surfaces
of genus at most that of the diagram $D$. Therefore, both have canonical genus at most 
$g$, and so have canonical genus $g$ (since $K$ does). But since the new diagrams 
for $K$ and $K_1$
have fewer crossings than $D$, after finitely many applications of this
process, the alternating knot $K^\prime$ corresponding to a diagram of $K$ is prime, 
hence hyperbolic, and Seifert's algorithm applied to the diagram $D$ giving both
$K$ and $K^\prime$ gives surfaces $\Sigma$ and $\Sigma^\prime$ of genus $g$.

\leavevmode

\epsfxsize=4in
\centerline{{\epsfbox{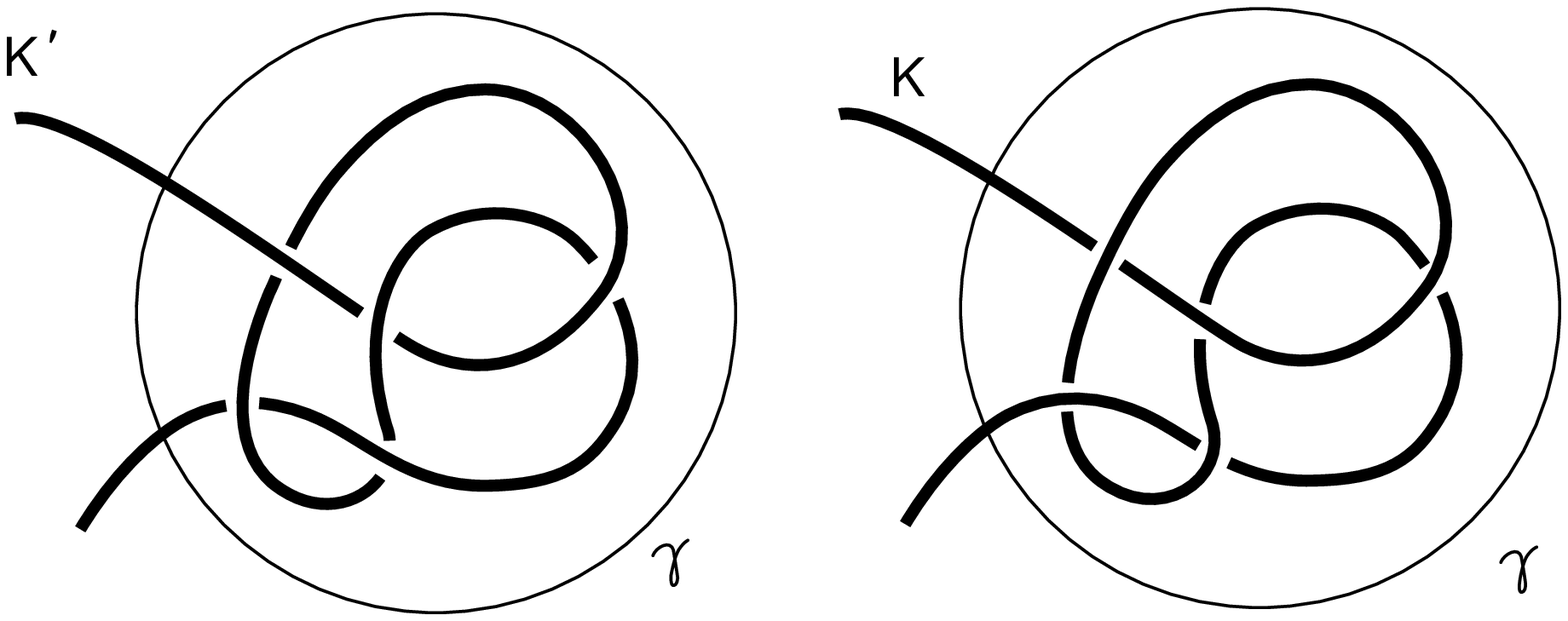}}}

\centerline{Figure 3}

Therefore, every hyperbolic knot $K$ with canonical genus $g$ has a diagram $D$ so 
that the corresponding alternating knot $K^\prime$, with diagram $D^\prime$, is 
hyperbolic and has genus, hence canonical genus, $g$.

\heading{\S 2 \\ Building the links $L_i$}\endheading

For each hyperbolic knot $K$ with canonical genus $g$, we have 
now realized 
this genus by a diagram $D$ so that the corresponding alternating
diagram $D^\prime$ represents a hyperbolic knot $K^\prime$. 
Now we will replace this knot $K^\prime$ by a hyperbolic link $L$ 
so that both $K$ and $K^\prime$ are realized as $1/n_i$ surgeries on 
the unknotted components of $L$. What we will see is that, as this
construction is carried out for all such knots $K$, we will in 
fact built only finitely many different links $L$.

The basic idea is that the Seifert circles for $\Sigma^\prime$ (and hence for 
$\Sigma$; they
are the same) come in two types; those that meet at most two half-twisted
bands, and those that meet three or more. The essential idea is to more or less
throw away the Seifert circles of the first type, and use the remainder to build
our link $L$. We can think of these Seifert surfaces for $K^\prime$ and $K$ as
being constructed from the (possibly non-planar) graph $\Gamma$ with fat vertices. 
The vertices of $\Gamma$
are the Seifert circles, and the edges of $\Gamma$ represent the half-twisted 
bands (see Figure 4). Such a graph is often called a rigid-vertex graph. The
Seifert circles of the first type therefore represent the vertices of 
valency one and two.

\leavevmode

\epsfxsize=4.9in
\centerline{{\epsfbox{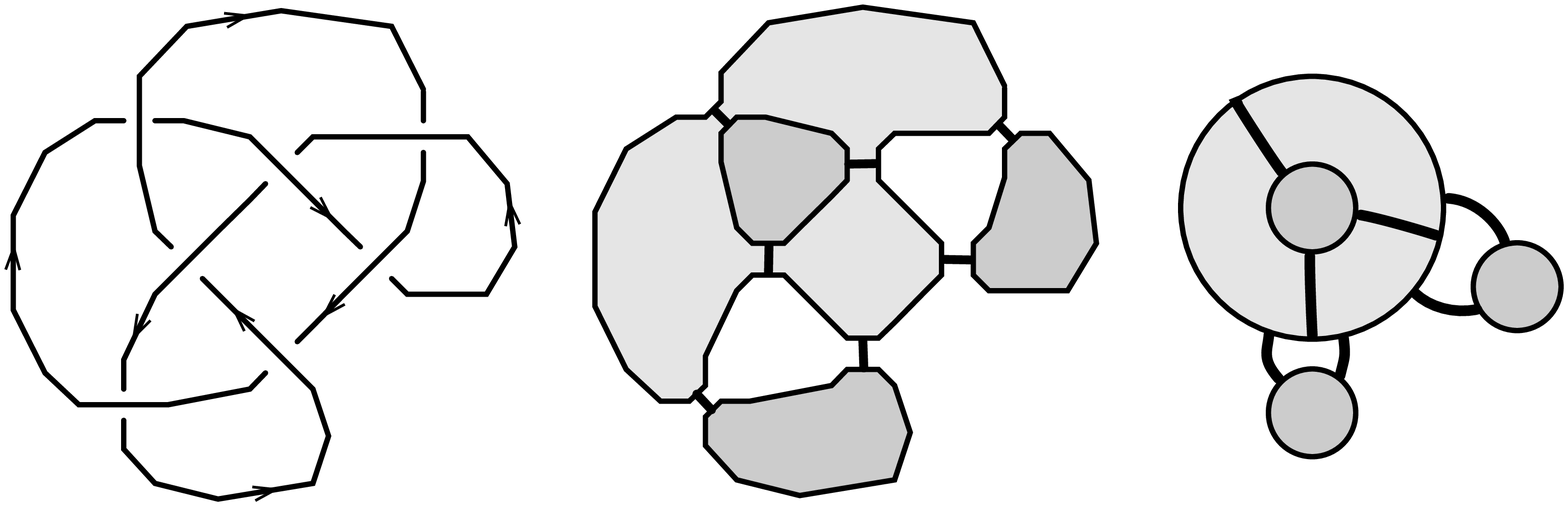}}}

\centerline{Figure 4}

The Seifert circles meeting only one half-twisted band can be safely ignored; 
they represent nugatory crossings of $K$. Removing the Seifert circle and the 
associated half-twisted band yields a new diagram for our knot $K$ with one fewer
crossing. Running Seifert's algorithm on it yields our original surface $\Sigma$ with
the disk and half-twisted band removed, which has the same genus as, and in fact is 
isotopic to, our original surface. The same is true of the alternating knot 
$K^\prime$.

Note that $\Sigma^\prime$ must have a Seifert circle meeting at least three bands; 
the only links whose (canonical) Seifert surfaces have no such circles 
are the (2,2$n$) torus links. 
The disks meeting two bands, together with the half-twisted bands that they meet, 
therefore form arcs in the the graph $\Gamma$. The nature of Seifert's algorithm
dictates that these arcs do not cross one another. The number of vertices in each
arc can be even or odd, depending on whether the normal orientations of the Seifert
circles at each end agree or differ, and whether or not the circles are nested. 
The half-twisted bands in 
each arc all twist in the same direction, because our knot $K^\prime$ is
alternating.

\leavevmode

\epsfxsize=3.5in
\centerline{{\epsfbox{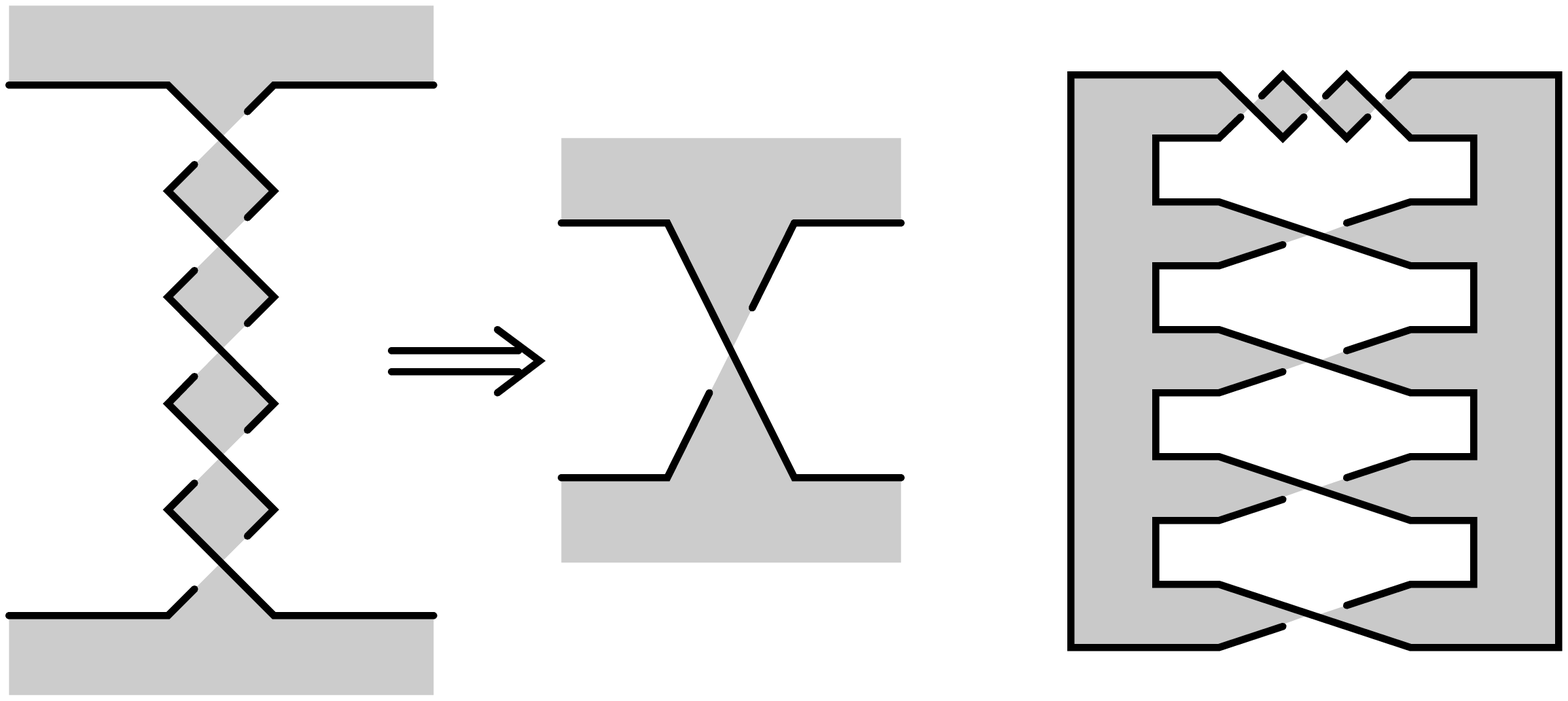}}}

\centerline{Figure 5}

What we will do now is replace each arc of 
circles and bands by a single circle with two twisted bands, if there were an
even number of bands in the arc, and a single twisted band, if there
were an odd number (Figure 5), each time 
having the bands twist in the same direction that we started with
in our alternating knot $K^\prime$. Also, for insurance, 
if all arcs were replaced by single twisted bands, 
choose one and replace it instead with two disks and three twisted bands. This is to avoid
inadvertently building a (2,n) torus knot; this would otherwise happen if we apply the
above procedure to a pretzel-type knot with all odd twists (see Figure 5), i.e., if there were 
only two Seifert circles meeting more than two bands, and every arc of bands ran from 
one to the other.

This new surface $\Sigma_0$ is the same surface that 
Seifert's algorithm would build for the underlying diagram, 
since we have respected the normal orientations  of the disks at the ends. 
It also has genus $g$, since all we have done is take out several full twists
in each arc of Seifert circles, so the surface is homeomorphic to the surface we started with.
The boundary $K_0$ of $\Sigma_0$ is also alternating, since we have kept the direction 
of twisting intact. This knot $K_0$ is therefore hyperbolic, since a loop $\gamma$ in 
the projection plane
meeting the diagram twice, as in section 1, would do the same for the 
diagram $D$, since 
we have only removed full twists; we simply imagine concentrating our added twists
in the vicinity of one of the crossings that remained, and no new intersections with
$\gamma$ will be created. Since in the original diagram $D$ this loop must cut off a 
trivial arc, it must do so for our new diagram, as well; the trivial arc has no
crossings that might have been removed, so it remains intact.

\smallskip

From this knot $K_0$ we build our link $L$ by placing a loop around each string of 
(one or two or three)
twisted bands, disjoint from the Seifert surface (Figure 6). This link is then an augmented 
alternating link, in the sense of Adams [Ad], and so is hyperbolic. Furthermore, since
$1/n$ surgery on one of the added unknotted components of $L$ amounts, according to the
Kirby calculus [Ro], to adding $n$ full twists to the arc of bands, we can recover
our alternating knot $K^\prime$ by doing such surgeries on the added components. On the
other hand, $K$ can be obtained from $K^\prime$ by changing crossings, which simply amounts
to adding another full twist to the half-twisted band representing the crossing (Figure 6).
Therefore, by adding up the net number of full twists along each string of bands,
our original knot $K$ can also be obtained by doing $1/n_i$ Dehn surgeries on the 
unknotted components of $L$.

\leavevmode

\epsfxsize=4in
\centerline{{\epsfbox{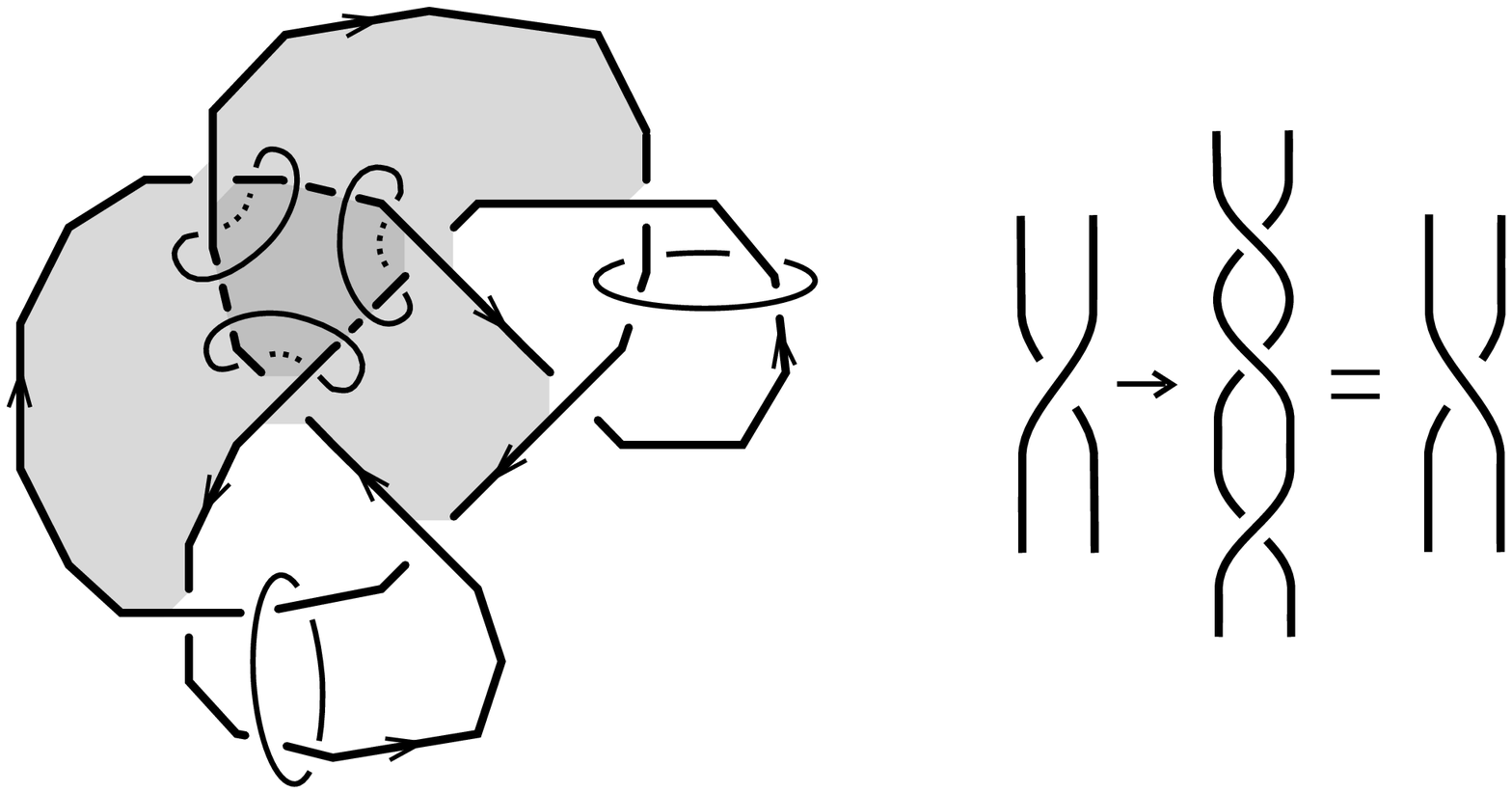}}}

\centerline{Figure 6}

\smallskip

Note that, in the end, our device of adding three twisted bands to one of the arcs is unnecessary.
The above argument shows that the corresponding link is hyperbolic, but the link complement
is homeomorphic to the one we would have obtained by adding only one twisted band, and so is
also hyperbolic. We could not prove this directly, however, using the results of [Ad].

\heading{\S 3 \\ Counting the $L$'s}\endheading

It remains only to show that the construction above will build only finitely-many
distinct links $L$. But this is simply a matter of counting the number of crossings
in the underlying knot $K_0$. The essential point is that the crossings of $K_0$
occur only in the arcs of half-twisted bands, and these are all very short. The genus of 
the Seifert surface $\Sigma_0$, however, can be computed from the number of Seifert circles 
meeting three
or more twisted bands (i.e., of the second type), and this genus is bounded. 
Therefore, there can't be very many circles of the second type. Therefore, there can't be very
many arcs of twisted bands (again, because the genus is bounded), so there can't be 
very many crossings. 

More specifically,
let $N_i$ be the number of circles of the second type meeting exactly $i$ bands, and
$N$ the number of circles of the second type. Then the number of arcs of bands, $A$, is equal to
one-half of $\sum iN_i$, and so

\smallskip

\noindent $\chi(\Sigma_0$) = $1-2g$ = $N-A$ = $\displaystyle N-{{1}\over{2}}\sum iN_i$ =
$\displaystyle \sum N_i(1-{{1}\over{2}}i)$

\hfill $\leq$ $\displaystyle \sum N_i(1-{{3}\over{2}})$ = 
$\displaystyle -{{1}\over{2}}N$

\smallskip

\noindent so $N\leq 4g-2$. Consequently, 

\centerline{$1-2g$ = $\displaystyle N-{{1}\over{2}}\sum iN_i$ $\leq$ $\displaystyle (4g-2)-{{1}\over{2}}\sum iN_i$, so}

\smallskip

\centerline{ $A$ = $\displaystyle{{1}\over{2}}\sum iN_i$ $\leq$ $6g-1$ .}

\smallskip

\noindent Consequently, since each additional unknotted loop around an arc of bands
adds four crossings to the diagram of $L$, 

\smallskip

\noindent\#(crossings of $L$) = 4$A$ + \#(crossings of $K_0$) $\leq$
$4A+(2A+1)$ 

\hfill = $6A+1$ $\leq$ $6(6g-1)+1$ $\leq$ $36g$

\smallskip

\noindent since each arc of bands contributes at most two crossings to $K_0$, except possibly for 
one which contributes three. Because there is a bound on the number of crossings in the link $L$, 
there can therefore be only finitely many such links, and therefore there is a finite bound $C(g)$
on their volume. Every hyperbolic knot with canonical genus $g$ is obtained by doing $1/n_i$ 
surgeries on the unknotted components of one of these links, and therefore has volume less than 
$C(g)$, as well. This completes the proof of the theorem.

\heading{\S 4 \\ A linear bound on volume}\endheading

The fact that we can bound the number of crossings in the link $L$ that we build $K$ from
gives us a way to give a linear bound on the volume of $K$ in terms of the genus $g$.
The author is indebted to William Thurston for pointing this out.
The starting point is the fact that we can use a diagram of a knot to ideally triangulate its
complement; see the documentation for SnapPea [We] for an outline of the process. The
resulting triangulation has at most 4 times as many tetrahedra as the diagram has crossings.
The idea is basically to triangulate the top and bottom hemispheres of the 3-sphere with 
vertices at the poles,
and with all other vertices lying on the knot projection, largely by coning off the 
cell decomposition
of the projection sphere given by the knot diagram $D$, adding additional edges to 
cut the cells into triangles. The ideal
triangulation is obtained by collapsing a pair of edges to bring the poles to the knot, which 
in the knot complement is at infinity.

By pulling the edges tight, making them geodesics, the complement of the knot is then covered
(not triangulated; this is the phenomenon of `negatively oriented tetrahedra' found in SnapPea)
by geodesic ideal tetrahedra, each of which has volume less than the volume 
$V_0$ of a regular
ideal tetrahedron, which is slightly larger than 1. The volume of the link $L$ is therefore
bounded by $4V_0$ times the number of crossings of $L$.

We can slightly improve our bound on the number of crossings above, by noting that at every place
where an arc of bands was replaced by two twisted bands, we can, since we have also added the loop
travelling around the bands, remove the two crossings that the bands created; the resulting link
complement homeomorphic to the original one (so is still hyperbolic), and all of our knots $K$
can be obtained from these new links; we simply add another full twist around the 
unknotted loop. This lowers the
number of crossings to something bounded by $5A+2$ (the 2 is for the arc we might have replaced by
three twisted bands), giving a bound, in terms of the genus, of $5(6g-1)+2$ $\leq$ $30g$. The 
resulting bound on volume is then $4V_0(30g)$ $\leq$ $122g$

\leavevmode

\epsfxsize=4in
\centerline{{\epsfbox{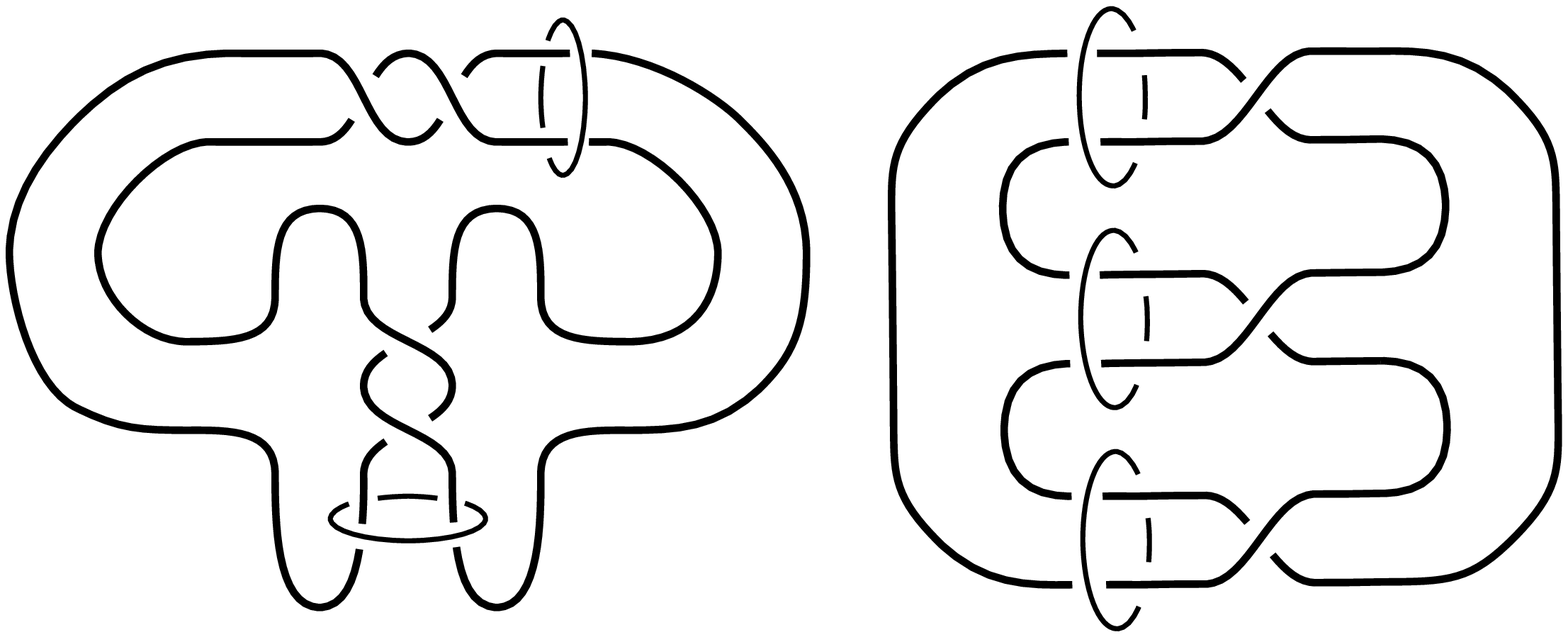}}}

\centerline{Figure 7}

\smallskip

This bound is of course very crude; for example, it is easy to see that the only knots with
canonical genus one are the 2-bridge knots of type (2$n$,2$m$) and the 
classical pretzel knots with
numbers of twists all odd; by Euler characteristic considerations, there can be at most two
Seifert circles which meet more than two twisted bands. 
These knots can be obtained by surgeries on the unknotted components of two links
(Figure 7), and so, according to SnapPea [We], the knots will all have volume less than 15.

\leavevmode

\epsfxsize=4in
\centerline{{\epsfbox{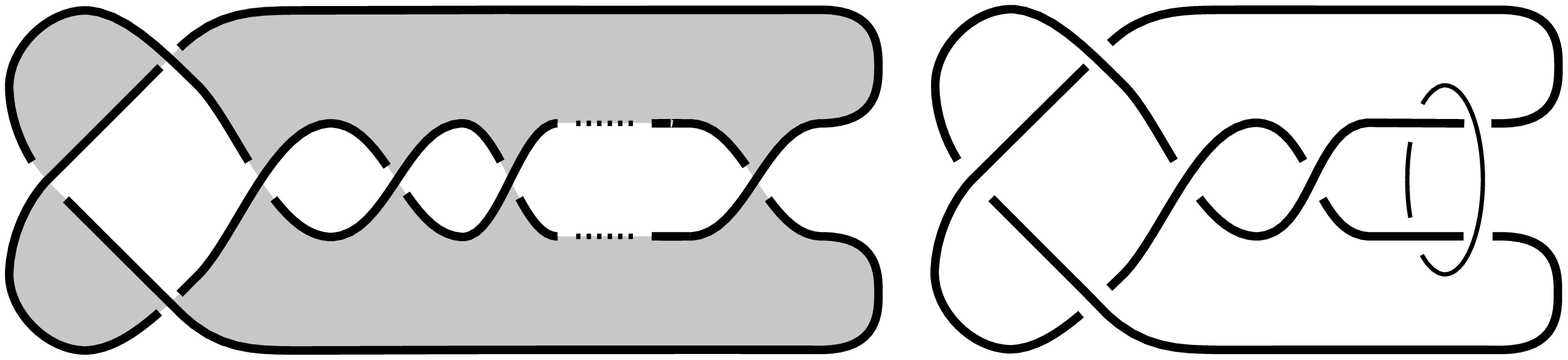}}}

\centerline{Figure 8}

\smallskip

We should note that it is possible for the genus, hence canonical genus, of a knot to be 
arbitrarily large, while its volume remains small. Figure 8 provides a collection of 
alternating knots whose canonical genera go to infinity. All of these knots,
however, can be 
obtained by $1/n$ Dehn surgeries on the unknotted component of the link in Figure 8, and
so their volumes are all bounded from above by the volume of the link.

The fact that we could find a volume bound at all in our theorem stemmed from the fact
that all of the knots in our family could be obtained by doing Dehn surgeries on a finite 
number of links. It seems natural to ask if there is a converse to this relationship.
Specifically, if we consider the family of hyperbolic knots whose volumes are bounded 
by a fixed constant, can they all be obtained by Dehn filling along components of some
finite number of links? (Formally, taking the union of these links, they would then
be obtained by Dehn fillings on the components of a single link.) This is perhaps not
entirely unreasonable, since the set of volumes is well-ordered, with order type 
$\omega^\omega$ [Th], and the function from knot complement to volume is finite to one. 
Using our understanding of how limit points arise in this set (as Dehn drilling), 
we might be able to identify some finite set of (limits of .... limits of limit) volumes, 
corresponding to a finite set of links, so that every knot
in our family shares its volume with a knot obtained by surgery on one of these links.
But this, of course, does not really suffice.
An added complication is the fact that we do not (necessarily) ask that the components
of the links that we do Dehn fillings along be unlinked from one another, as they
were in the construction above. The above approach, for example, would provide no
information about the linking of the added components.

\smallskip

\centerline{\bf Acknowledgements}

\smallskip

The author wishes to thank Jeff Weeks for many helpful conversations during the
course of this work, and Tsuyoshi Kobayashi and Nara Women's University for 
providing the setting in which these conversations took place.

\enddocument